\input amstex
\input Amstex-document.sty


\define\Vol{\operatornamewithlimits{Vol}}
\define\conv{\operatornamewithlimits{conv}}
\define\Area{\operatornamewithlimits{Area}}

\define\Prob{\operatornamewithlimits{Prob}}

\pageno 527

\def\shorttitle{Random Points, Convex Bodies, Lattices}
\def\shortauthor{I. B\'ar\'any}

\topmatter

\title\nofrills{\boldHuge Random Points, Convex Bodies, Lattices}
\endtitle

\author {\Large Imre B\'ar\'any*} \endauthor

\thanks {*R\'enyi Institute of Mathematics, Hungarian Academy of Sciences,
PoB 127, Budapest 1364, Hungary, and Department of Mathematics, University College London, Gower Street, London
WC1E 6BT, UK. E-mail: barany\@renyi.hu and barany\@math.ucl.ac.uk}\endthanks

\abstract\nofrills \centerline{\boldnormal Abstract}

\vskip 4.5mm

{\ninepoint Assume $K$ is a convex body in $R^d$, and $X$ is a (large) finite subset of $K$. How many convex
polytopes are there whose vertices come from $X$? What is the typical shape of such a polytope? How well the
largest such polytope (which is actually $\conv X$) approximates $K$? We are interested in these questions mainly
in two cases. The first is when $X$ is a random sample of $n$ uniform, independent points from $K$ and is
motivated by Sylvester's four-point problem, and by the theory of random polytopes. The second case is when $X=K
\cap Z^d$ where $Z^d$ is the lattice of integer points in $R^d$. Motivation comes from integer programming and
geometry of numbers. The two cases behave quite similarly.

\vskip 4.5mm


}
\endabstract
\endtopmatter

\document

\baselineskip 4.5mm \parindent 8mm

\specialhead \noindent \boldLARGE 1. Sylvester's question \endspecialhead

In  the 1864 April  issue of  the  Educational Times  J. J.
Sylvester [26] posed the innocent looking  question  that read:
{\it ``Show   that the  chance  of four points forming the apices
of a reentrant quadrilateral is 1/4 if they be taken at random in
an indefinite plane.''} It was understood within  a year that the
question is ill-posed. (The culprit  is, as  we all  know by now,
the ``indefinite plane'' without a properly defined probability
measure on it.) So Sylvester modified the  question: let $K\subset
R^2$ be a convex body (that is, a compact, convex set with
nonempty interior) and choose four random, independent points
uniformly from $K$, and write $P(K)$ for the probability that the
four points form the apices of a reentrant  quadrilateral, or,  in
more modern  terminology, that their convex hull is a  triangle.
How large is  $P(K)$,  and for what  $K$  is $P(K)$ the largest
and the smallest? This question became known as Sylvester's
four-point problem. It took fifty years to find the  answer:
Blaschke [16] showed that for all convex bodies $K \subset R^2$
$$
P(\text{disk})\leq P(K)\leq P(\text{triangle}).
$$

Assume now, more generally, that $X_n=\{x_1,\dots,x_n\}$ is a random sample of $n$ uniform,  independent points
from the  convex body $K$  and write $p(n,K)$ for the probability that $X_n$ is in {\sl convex position}, that is,
no $x_i$ is in the    convex  hull of the   others.   Sylvester's question   is  just  the complementary   problem
for  $n=4$:   $P(K)=1-p(4,K)$. The probability $p(n,K)$ has been determined in various special cases (see [22, 17,
9, 27]). The following  result from [7] describes the asymptotic behaviour of $p(n,K)$.

{\bf Theorem 1.1.} \it For every convex body $K \subset R^2$ of
unit area
$$
\lim_{n \to \infty} n^2 \root n \of {p(n,K)}= \frac {e^2}{4}
A^3(K)$$ where $A(K)$ is the supremum of the affine perimeter of
all convex sets $S \subset K$. \rm

The affine perimeter, $AP(K)$ can be defined in many ways (see
[23]), for instance $AP(K)=\int _{\partial K} \kappa^{1/3}ds$
where $\kappa$ is the curvature and integration goes by
arc-length. (This definition works for smooth convex bodies, the
extension for all convex bodies can be found in [23].) Theorem 1
of [6] says that there is a unique convex compact set $K_0 \subset
K$ with $AP(K_0)=A(K)$. The proof of Theorem 1.1 gives more than
just the asymptotic behaviour of $p(n,K)$, namely, if the random
points $x_1,\dots,x_n$ are in convex position, then their convex
hull is, with high probability, very close to $K_0$. For the
precise formulation see [7].

Define $Q(X_n)$ as the collection of all convex polygons spanned
by the points of $X_n$, that is, $P \in Q(X_n)$ iff
$P=\conv\{x_{i_1},\dots,x_{i_k}\}$ for some $k$-tuple of points
from $X_n$ that is in convex position ($k\geq3$). Clearly,
$Q(X_n)$ is a random collection as it depends on the random sample
$X_n$. How many polygons are there in $Q(X_n)$? The answer is
given in [7]. Write $E|Q(X_n)|$ for the expectation of the size of
$Q(X_n)$.

{\bf Theorem 1.2.} \it For every convex body $K \subset R^2$ of
unit area
$$
\lim_{n \to \infty} n^{-1/3} \log E|Q(X_n)|= 3\cdot 2^{-2/3}A(K).
$$
\rm

Further, there is a limit shape to the polygons in $Q(X_n)$,
meaning that all but a small fraction of the polygons in $Q(X_n)$
are very close to $K_0$. We use $\delta (S,T)$ to denote the
Haussdorf distance of $S,T \subset R^2$.

{\bf Theorem 1.3.} \it For every convex body $K \subset R^2$ and
for every $\varepsilon >0$
$$
\lim_{n \to \infty} \frac {E|\{P \in Q(X_n)\: \delta
(P,K_0)>\varepsilon\}|}{E|Q(X_n)|}=0.
$$ \rm

We will see in Section 3 that similar phenomena hold for the
lattice case. In general, lattice points and random points, in
relation to convex bodies, behave very much alike. Quite often one
understands in the random case what to expect for lattice points,
or the other way around. The proofs are quite different and are
omitted in this survey.

\specialhead \noindent \boldLARGE 2. Higher dimensions
\endspecialhead

Much less is known in higher dimensions. One reason is that the
unicity of the convex subset of $K$ with maximal affine surface
area is not known. It is a mystery, for instance, which convex
subset of the unit cube in $R^3$ has maximal affine surface area.
But there are other reasons as well, connected to the lack of the
multiplicative rule (5.3) from [7]. Yet one can prove the
following asymptotic formula [8]. Here $\Cal C ^d$ denotes the set
of all convex bodies in $R^d$, and $p(n,K)$ denotes, as before,
the probability that the random sample $X_n=\{x_1,\dots,x_n\}$
from $K$ is in convex position, that is, no $x_i$ is in the convex
hull of the others.

{\bf Theorem 2.1.} \it For every $K \in \Cal C ^d$, and for all
$n\geq n_0$
$$
c_1 < n^{2/(d-1)} \root n \of {p(n,K)} <c_2
$$
where $n_0,c_1,c_2$ are positive constants that depend only on $d$. \rm

With Vinogradov's convenient $\ll$ or $\ll_d$ notation this says
that
$$
1 \ll_d n^{2/(d-1)} \root n \of {p(n,K)} \ll_d 1.
$$
From this one can estimate the size of $E|Q(X_n)|$ when $K \in
\Cal C^d$:
$$
n^{(d-1)/(d+1)} \ll_d \log E|Q(X_n)| \ll_d n^{(d-1)/(d+1)}.
$$

For comparison let us have a look at lattice polytopes contained
in some fixed $K \subset \Cal C ^d$. So let $Z^d$ be the lattice
of the integers in $R^d$ and consider, for a large integer $m$,
the lattice $\frac 1m Z^d$. Assume $K$ contains $n$ points from
this lattice. As $m$ is large, $n=(1+o(1))m^d \Vol K$. Write $\Cal
P_m(K)$ for the collection of all $\frac 1m Z^d$-lattice polytopes
contained in $K$. The next theorem, which follows easily from the
results of [12], shows a very strong analogy between $\Cal P_m(K)$
and $Q(X_n)$.

{\bf Theorem 2.2.} \it For every $K \in \Cal C ^d$
$$
n^{(d-1)/(d+1)}\ll_d \log |\Cal P _m(K)| \ll_d n^{(d-1)/(d+1)}.
$$
\rm

The result shows that when $K \in \Cal C^d$ contains $n$ lattice
points, these lattice points span (essentially) $\exp
\{cn^{(d-1)/(d+1)}\}$ convex polytopes, the same number as in the
random case. Lattice points and random points in convex bodies
behave similarly: this is the moral.

\specialhead \noindent \boldLARGE 3. Lattice polygons and limit shape \endspecialhead

In the plane Theorem 2.2 can be proved in stronger form (see [5],
[6], and [28]):

{\bf Theorem 3.1.} \it For every $K \in \Cal C^2$
$$
\lim n^{-2/3} \log |\Cal P_n(K)|=3\root 3 \of {\frac
{\zeta(3)}{4\zeta(2)}}A(K).
$$\rm

Here $\zeta(.)$ stands for Riemann $\zeta$ function. Note that
this result is in complete analogy with Theorem 1.2: just the
constant is different. (Also, $n$ is in power $-2/3$ instead of
$-1/3$ as $K$ contains $(1+o(1))n^2\Area K$ lattice points.) The
analogy carries over to Theorem 1.3 as well:

{\bf Theorem 3.2.} \it  For every convex body $K \in \Cal C^2$ and
for every $\varepsilon >0$
$$
\lim \frac {|\{P \in \Cal P_n(K))\: \delta
(P,K_0)>\varepsilon\}|}{|\Cal P_n(K)|}=0.
$$ \rm

This shows again that all but a tiny fraction of the polygons in
$\Cal P_n(K)$ are very close to $K_0$. In other words, these
polygons have a limit shape. Theorems of this type were first
proved by B\'ar\'any [5], Vershik [28] (for the case when $K$ is
the unit square). Sinai [25] found a different proof which uses
probability theory and gives a central limit theorems about how
small that tiny fraction of polygons is. This has been generalized
by Vershik and Zeitouni [29] to all convex bodies in $R^2$. A
central limit theorem of this type holds for the random sample
case as well, see [14] for the precise statement.

\specialhead \noindent \boldLARGE 4. The integer convex hull \endspecialhead

The {\sl integer convex hull}, $I(K)$, of a convex body $K \in
\Cal C^d$ is, by definition, the convex hull of the lattice points
contained in $K$:
$$
I(K)=\conv(Z^d \cap K).
$$
$I(K)$ is clearly a convex polytope. How many vertices does it
have? Motivation for the question comes from integer programming,
classical enumeration questions (like the circle problem), and
from the theory of random polytopes. In integer programming one
wants to know that $I(K)$ does not have too many vertices,
assuming, say, that $K$ is a nice rational polytope. The latter
means that $K$ can be given by $m$ inequalities with integral
coefficients; the size of such an inequality is the number of bits
necessary to encode it as a binary string. Then the size of the
rational polytope is the sum of the sizes of the defining
inequalities. Strengthening earlier results by Shevchenko [24],
and Hayes and Larman [21], Cook, Hartman, Kannan, and McDiarmid
[18] showed that for a rational polytope $K$ of size $\phi$
$$
f_0(I(K)) \leq 2m^d(12d^2\phi)^{d-1}.
$$
Here, as usual, $f_i(P)$ stands for the number of $i$-dimensional
faces of the polytope $P$. Thus $f_0(P)$ is the number of vertices
of $P$. Most likely, the same inequality holds for all
$i=0,1,\dots ,d-1$:
$$
f_i(I(K)) \ll \phi^{d-1}
$$
where the implied constant depends on $d$ and $m$ as well.

The above inequality for $f_0(I(K))$ is best possible, as is shown
B\'ar\'any, Howe, and Lov\'asz in [11]:

{\bf Theorem 4.1.} \it For fixed $d\geq 2$ and for every $\phi >0$ there exists a rational simplex $P \subset R^d$
of size at most $\phi$ such that $I(P)$ has $\gg_d \phi^{d-1}$ vertices. \rm

The construction uses algebraic number theory. It shows further
that the estimate $f_i(I(K)) \ll \phi^{d-1}$ for all $i$
 is best possible, if true.

What about the integer convex hull of other convex bodies? Balog
and B\'ar\'any [3] considered case $K=rB^2$ where $B^2$ is the
Euclidean unit ball centered at the origin and $r$ is large and
showed that
$$
0.3r^{2/3}< f_0(I(rB^2))<5.5r^{2/3}.
$$
Later Balog and Deshoullier [4] determined the average of
$f_0(I(rB^2))$ on an interval $[R,R+H]$ which turned out to be
very close to $3.453R^{2/3}$ as $R$ goes to infinity ($H$ has to
be large). B\'ar\'any and Larman [13] determined the order of
magnitude of $f_i(I(rB^d))$. (The method, and the result, apply
not only to the unit ball but to smooth enough convex bodies as
well.)

{\bf Theorem 4.2.} \it For every $d \geq 2$ and every
$i=0,1,\dots,d-1$
$$
r^{d(d-1)/(d+1)} \ll _d f_i(I(rB^d)) \ll_d r^{d(d-1)/(d+1)}.
$$ \rm

This result is related to a beautiful theorem of G. E. Andrews [1]
stating that a lattice polytope $P$ in $R^d$ with volume $V>0$ has
$\ll_d V^{(d-1)/(d+1)}$ vertices. The above theorem shows that
Andrews' estimate is best possible (apart from the constant
implied by $\ll _d$). A similar (perhaps less compact) example was
given earlier V. I. Arnol'd [2].

This kind of question can be considered in a more general setting.
Let $G$ be the group of all isometries of $R^d$ with translations
by elements of $Z^d$ factored out. $G$ is a compact topological
group with a Haar measure which is a unique invariant probability
measure when normalized properly. Assume $g \in G$ is chosen
according to this probability measure. Then $gK$ is a random copy
of $K$ and we can talk about the expectation of the random
variable $f_0(I(gK))$.

For the next result we assume $K \in \Cal C^d$ and define the
function $u\:K \to R$ by
$$
u(x)= \Vol (K \cap (2x-K)),
$$
that is, $u(x)$ is the volume of the intersection of $K$ with $K$
reflected about $x$. Set, finally, $K(u<t)=\{x \in K \: u(x)<t\}$.
The following is an unpublished result of B\'ar\'any and Matou{\v
s}ek:

{\bf Theorem 4.3.} \it Consider all $K \in \Cal C ^d$ with the
ratio of the radii of the smallest circumscribed and the largest
inscribed balls to $K$ bounded by $D$. Then, as $\Vol K$ goes to
infinity,
$$
\Vol K(u<1) \ll Ef_0(I(gK)) \ll \Vol K(u<1)
$$
where the constants implied by $\ll$ depend only on $d$ and $D$. \rm

It follows easily from Minkowski's classical theorem that all
vertices of $I(K)$ belong to $K(u<2^d)$. (This is the first step
in proving the upper bound.) It is not hard to see that $\Vol
K(u<2^d) \ll \Vol K(u<1)$. So the meaning of the theorem is that
the average number of vertices of $I(gK)$ is essentially the
volume of $K(u<2^d)$. Probably the same is true for the expected
number of $i$-dimensional faces of $I(gK)$ but there is no proof
in sight.

The behaviour of $\Vol K(u<1)$ is more or less known (from [10],
say, but more precise results are known as well): it is of order
$(\Vol K)^{(d-1)/(d+1)}$ for smooth enough convex bodies and of
order $(\log \Vol K)^{d-1}$ for polytopes, and it is between these
bounds for all convex bodies.

We mention further that Theorem 4.3 is quite analogous to a result
in the theory of random polytopes. Given $K \in \Cal C ^d$, and a
random sample of $n$ points, $X_n$, from $K$, $K_n= \conv X_n$ is
called a {\sl random polytope} on $n$ points. It is shown in [10]
that, assuming $\Vol K=n$ (which is the proper scaling for
comparison with Theorem 4.3), for all $i=0,1,\dots,d-1$
$$
\Vol K(u< 1) \ll Ef_i(K_n) \ll \Vol K(u < 1)
$$
where the implied constants depend only on dimension.

Note that, unlike Theorem 4.3, this result works for all
$i=0,\dots, d-1$ (without any condition on the ratio of radii of
the circumscribed and inscribed balls). Most likely, Theorem 4.3
also holds for all $i$, which would make the analogy even more
complete.

There is, however, a point here where the analogy breaks down. Let
$K\subset R^2$ be the square of area $n$, so $K_n$ is a random
polytope, and $I(gK)$ is the integer hull of a random copy of $K$.
The expectation of $\Area (K \setminus K_n)$ is of order $\log n$
(see [10], say), while the expectation of $\Area (K \setminus
I(gK))$ is of order $(\log n)^2$. (The latter result comes again
from the unpublished work of B\'ar\'any and Matou{\v s}ek.) The
reason is that the boundary of $K_n$ contains no points from $X_n$
apart from its vertices, while the boundary of $I(gK)$ does. A
further reason is that what we are measuring here is a metric
property, and not a combinatorial one. We think that the same
phenomena is bound to happen in higher dimension.

\specialhead \noindent \boldLARGE 5. Random 0-1 polytopes \endspecialhead

Finally we mention a recent development, prompted by a question of
K. Fukuda and G. M. Ziegler [30]. They asked how many facets a 0-1
polytope in $R^d$ can have; a 0-1 polytope is a polytope whose
vertices only have 0 or 1 coordinates. So such a polytope is the
convex hull of a subset of the vertices of the unit cube, $Q^d$,
in $R^d$. 0-1 polytopes play an important role in combinatorial
optimization where the target is, very often, a concise
description of the facets of the polytope. This task has turned
out to be difficult for several classes of 0-1 polytopes.

Write $G(d)$ for the maximal number of facets a 0-1 polytope can
have. It is not hard to see that $2^d \leq G(d) \leq 2d!$. The
upper and lower bounds have been improved slightly: the lower
bound by a construction of Christoff (see [30]), and the upper
bound by Fleiner, Kaibel, and Rote [20].

The vertices of every 0-1 polytope are on a sphere (centered at
$(1/2,\dots,1/2)$). There is a formula (see for instance [9]) for
the expected number of facets of a random polytope with $n$
uniform independent points from the (unit) sphere in $R^d$. It
says that, in the range when $2d <n<2^d$, the expected number of
facets is of order $(\log n/d)^{d/2}$. So if the analogy between
random points and lattice points carries over the 0-1 case one
should expect $G(d)$ to be of order $d^{d/2}$. This is too much to
ask for at the moment, yet the following is true (see [15]).

\bf Theorem 5.1. \it There is a constant $c>0$ such that for all $d\geq 2$
$$
G(d) \gg \left(c \frac {d}{\log d}\right)^{d/4}.
$$ \rm

The construction giving this estimate is random. Write $K_n$ for
the convex hull of $n$ random, uniform, and independent 0-1
vectors. Assume $x$ is a point from $Q^d$, and define
$$
p(x,n) = \Prob[x \in K_n ].$$
 General principles would tell that,
for most $x \in Q^d,\; p(x,n)$ is either close to one or close to
zero. To be more specific, set
$$
P(t) = \{ x \in Q^d : p(x,n) \geq t \}.
$$
The proof of Theorem 5.1 is based on the fact that for all small
$\varepsilon> 0$ and large enough $d$  $P(1-\varepsilon ) \subset
P(\varepsilon  )$, of course, but the drop from $1-\varepsilon$ to
$\varepsilon$ is very abrupt: $P(\varepsilon )$ is in a small
neighbourhood of $P(1-\varepsilon )$. This shows that
$P(1-\varepsilon ) \subset K_n$ with high probability. But only a
tiny fraction of $K_n$ lies outside $P(\varepsilon )$ : most of
the boundary of $P(\varepsilon )$ is outside $K_n$. Thus most of
the boundary of $P(\varepsilon )$ is cut off by facets of $K_n$.
These facets lie outside $P(1-\varepsilon )$. Comparing the
surface area of $P(\varepsilon )$  with the amount a facet can cut
off from it gives the lower bound.

The actual proof is technical, difficult, and makes extensive use
a beautiful result of Dyer, F\"uredi, and McDiarmid [19]. Their
target was to determine the threshold $n=n(d)$ such that $K_n$
contains most of the volume of $Q^d$. As they prove, this happens
at $n= (2/\sqrt e )^d$. Their method describes where $p(x,n)$
drops from one to zero as $d \rightarrow \infty$. The analysis
carries over for other values of $n$. In our case higher precision
is required as we need a good estimate on how fast $p(x,n)$ drops
from one to zero. We were able to control this only where  the
curvature of the boundary of $P(\varepsilon )$ behaves nicely.
This is perhaps the spot where the exponent $d/2$ (for the random
spherical polytope) is lost and we only get $d/4$ for $K_n$.

\widestnumber \key {[99]}

\specialhead \noindent \boldLARGE References \endspecialhead

\ref \key 1 \by G. E. Andrews \paper \rm A lower bound for the
volume of strictly convex bodies with many boundary points \jour
\it Trans. Amer. Math. Society, \vol \rm 106 \yr 1963 \pages
270--279
\endref

\ref \key 2 \by V. I. Arnol'd \paper \rm Statistics of integral
lattice polytopes (in Russian) \jour \it Funk. Anal. Pril, \vol
\rm 14 \yr 1980 \pages 1--3
\endref

\ref \key 3 \by A. Balog, I. B\'ar\'any \paper \rm On the convex
hull of the integer points in a disc \jour \it Discrete Comp.
Geometry, \vol \rm 6 \yr 1992 \pages 39--44
\endref

\ref \key 4 \by A. Balog, J-M. Deshouliers \book \rm On some
convex lattice polytopes, in: {\it Number theory in progress,} \rm
(K. Gy\H ory, ed.)  \publ de Gruyter \yr 1999, 591--606
\endref

\ref \key 5 \by I. B\'ar\'any \paper \rm The limit shape of convex
lattice polygons \jour \it Discrete Comp. Geometry, \vol \rm 13
\yr 1995 \pages 270--295
\endref

\ref \key 6 \by I. B\'ar\'any \paper \rm  Affine perimeter and
limit shape \jour \it J. Reine Ang. Mathematik, \vol \rm 484 \yr
1997 \pages 71--84
\endref

\ref \key 7 \by I. B\'ar\'any \paper \rm Sylvester's question: the
probability that $n$ points are in convex position \jour \it
Annals of Probability, \vol \rm 27 \yr 1999 \pages 2020--2034
\endref

\ref \key 8 \by I. B\'ar\'any \paper \rm A note on Sylvester's
four point problem \jour \it Studia Math. Hungarica, \vol \rm 38
\yr 2001 \pages 73--77
\endref

\ref \key 9 \by I. B\'ar\'any, Z. F\"uredi \paper \rm On the shape
of the convex hull of random points \jour \it Prob. Theory Rel.
Fields, \vol \rm 77 \yr 1988 \pages 231--240
\endref

\ref \key 10 \by  I. B\'ar\'any, D. G. Larman \paper \rm Convex
bodies, economic cap coverings, random polytopes \jour \it
Mathematika, \vol \rm 35 \yr 1988 \pages 274--291
\endref

\ref \key 11 \by I. B\'ar\'any, R. Howe. L. Lov\'asz \paper \rm On
integer points in polyhedra: a lower bound \jour \it
Combinatorica, \vol \rm 12 \yr 1992 \pages 135--142
\endref

\ref \key 12 \by I. B\'ar\'any, A. M. Vershik \paper \rm  On the
number of convex lattice polytopes \jour \it GAFA Journal, \vol
\rm 2 \yr 1992 \pages 381--393
\endref

\ref \key 13 \by I. B\'ar\'any, D. G. Larman \paper \rm The convex
hull of the integer points in a large ball \jour \it Math.
Annalen, \vol \rm 312 \yr 1998 \pages 167--181
\endref

\ref \key 14 \by I. B\'ar\'any, G. Rote, W. Steiger, C. Zhang
\paper \rm A central limit theorem for random convex chains \jour
\it Discrete Comp. Geometry, \vol \rm 23 \yr 2000 \pages 35--50
\endref

\ref \key 15 \by I. B\'ar\'any, A. P\'or \paper \rm On 0-1
polytopes with many facets \jour \it Advances in Math., \vol \rm
\yr 2000 \pages 1--28
\endref

\ref \key 16 \by W. Blaschke \book \it Vorlesungen \"uber
Differenzialgeometrie II. Affine Differenzialgeometrie \publ
Springer \yr 1923
\endref

\ref \key 17 \by C. Buchta \paper \rm On a conjecture of R.E.
Miles about the convex hull of random points \jour \it Monatsh.
Math., \vol \rm 102 \yr 1986 \pages 91--102
\endref

\ref \key 18 \by W. Cook, M. Hartman, R. Kannan, C. McDiarmid
\paper \rm On integer points in polyhedra \jour \it Combinatorica,
\vol \rm 12 \yr 1992 \pages 27--37
\endref

\ref \key 19 \by  M. E. Dyer, Z. F\"uredi, C. McDiarmid \paper \rm
Volumes spanned by random points in the hypercube \jour \it Random
Structures and Algorithms,  \vol \rm 3 \yr 1992 \pages 91--106
\endref

\ref \key 20 \by  T. Fleiner, V. Kaibel, G. Rote: \paper \rm Upper
bounds on the maximal number of facets of 0/1-polytopes \jour \it
European J. Combinatorics, \vol\rm 21 \yr 2000 \pages 121--130
\endref

\ref \key 21 \by A. C. Hayes and D. G. Larman \paper \rm The
vertices of the knapsack polytope \jour \it Discrete App. Math.,
\vol \rm 6 \yr 1983 \pages 135--138
\endref

\ref \key 22 \by B. Hostinsky \paper \rm Sur les probabilit\'es
g\'eom\'etriques \jour \it Publ. Fac. Sci. Univ. Brno, \yr 1925
\endref

\ref \key 23 \by  E. Lutwak \paper \rm Extended affine surface
area \jour \it Adv. Math., \vol \rm 85 \yr 1991 \pages 39--68
\endref

\ref \key 24 \by V. N. Shevchenko \paper \rm On the number of
extreme points in linear programming (in Russian) \jour \it
Kibernetika, \vol \rm 2 \yr 1981 \pages 133--134
\endref

\ref \key 25 \by Ya. G. Sinai \paper \rm Probabilistic approach to
analyze the statistics of convex polygonal curves (in Russian)
\jour \it Funk. Anal. Pril., \vol \rm 28 \yr 1994 \pages 41--48
\endref

\ref \key 26 \by J. J. Sylvester \paper \rm Problem 1491 \jour \it
The Educational Times, (London) \yr April 1864 \pages 1--28
\endref

\ref \key 27 \by P. Valtr \paper \rm The probability that $n$
random points in a triangle are in convex position \jour \it
Combinatorica, \vol \rm 16 \yr 1996 \pages 567--574
\endref

\ref \key 28 \by A. M. Vershik \paper \rm The limit shape for
convex lattice polygons and related topics (in Russian) \jour \it
Funk. Anal. Appl., \vol \rm 28 \yr 1994 \pages 16--25
\endref

\ref \key 29 \by A. M. Vershik, O. Zeitouni \paper \rm Large
deviations in the geometry of convex lattice polygons \jour \it
Israel J. Math., \vol \rm 109 \yr 1999 \pages 13--28
\endref

\ref \key 30 \by G. M. Ziegler \book \rm  Lectures on 0/1
polytopes, in: {\it Polytopes---Combinatorics and Computation} \rm
(G. Kalai and G. M. Ziegler, eds.), DMV-Seminars, \publ
Birkh\"auser-Verlag \yr 2000, 1--44
\endref

\enddocument